\documentclass[12pt,leqno,draft]{article}

\newtheorem{theorem}{Theorem}
\newtheorem{lemma}[theorem]{Lemma}
\newtheorem{proposition}[theorem]{Proposition}
\newtheorem{definition}[theorem]{Definition}
\newtheorem{corollary}[theorem]{Corollary}

\newcommand{\begintheorem}{\addtocounter{equation}{1}\begin{theorem}}
\newcommand{\beginlemma}{\addtocounter{equation}{1}\begin{lemma}}
\newcommand{\beginproposition}{\addtocounter{equation}{1}\begin{proposition}}
\newcommand{\begindefinition}{\addtocounter{equation}{1}\begin{definition}}
\newcommand{\begincorollary}{\addtocounter{equation}{1}\begin{corollary}}

\begin{document}

\title{Notes on normed algebras, 5}

\author{Stephen William Semmes	\\
	Rice University		\\
	Houston, Texas}

\date{}

\maketitle

	Let $\mathcal{A}$ be a finite-dimensional commutative algebra
over the complex numbers with nonzero multiplicative identity element
$e$.  Let $\mathcal{M}$ denote the set of homomorphisms from
$\mathcal{A}$ onto the complex numbers.  Thus we get a homomorphism
from $\mathcal{A}$ into the algebra of complex-valued functions on
$\mathcal{M}$, with respect to pointwise addition and multiplication
of functions, by sending an element $x$ of $\mathcal{A}$ to the
function on $\mathcal{M}$ whose value at a homomorphism $\phi :
\mathcal{A} \to {\bf C}$ is $\phi(x)$.

	Suppose that $\phi$, $\psi$ are distinct elements of
$\mathcal{M}$.  This means that there is an $x \in \mathcal{A}$ such
that $\phi(x) \ne \psi(x)$.  By choosing $y$ to be a linear
combination of $x$ and $e$ we may arrange it so that $\phi(y) = 1$ and
$\psi(y) = 0$.  If $E$ is a finite subset of $\mathcal{M}$, then we
can use this construction to show that for any fixed $\phi \in E$
there is a $z \in \mathcal{A}$ such that $\phi(z) = 1$ and $\psi(z) =
0$ for all $\psi \in E$ such that $\psi \ne \phi$.

	Thus, for any finite subset $E$ of $\mathcal{M}$ and any
complex-valued function $f$ on $E$ there is a $w \in \mathcal{A}$ such
that $\phi(w) = f(\phi)$ for all $\phi \in E$.  In particular the
number of elements of $E$ is at most the dimension of $\mathcal{A}$.
It follows that the number of elements of $\mathcal{M}$ is at most the
dimension of $\mathcal{A}$, and that our homomorphism from $\mathcal{A}$
into the algebra of functions on $\mathcal{M}$ has image equal to
the whole algebra of functions on $\mathcal{M}$.

	Suppose that $x \in \mathcal{A}$ is in the kernel of the
canonical homomorphism from $\mathcal{A}$ onto the algebra of
functions on $\mathcal{M}$.  This is the same as saying that $\phi(x)
= 0$ for all homomorphisms $\phi$ from $\mathcal{A}$ onto the complex
numbers.  Equivalently, the spectrum of $x$ consists of exactly $0$.
By the Cayley--Hamilton theorem, this is equivalent to saying that $x$
is nilpotent, which is to say that $x^n = 0$ for some positive integer
$n$.

	If $\|\cdot \|$ is a norm on $\mathcal{A}$ such that
$(\mathcal{A}, \|\cdot \|)$ is a normed algebra, then
\begin{equation}
\label{|phi(x)| le |x|}
	|\phi(x)| \le \|x\|
\end{equation}
for all $x \in \mathcal{A}$ and $\phi \in \mathcal{M}$.  Of course a
natural norm on the algebra of functions on $\mathcal{M}$ is defined
simply by saying that the norm of a function $f$ is equal to the
maximum of $|f(\phi)|$ over all $\phi \in \mathcal{M}$.  The algebra
of functions on $\mathcal{M}$ equipped with this norm is a normed
algebra.  The previous assertion (\ref{|phi(x)| le |x|}) says exactly
that the canonical homomorphism from $\mathcal{A}$ onto the algebra of
functions on $\mathcal{M}$ has norm less than or equal to $1$, which
is to say that the norm of the image of an element of $\mathcal{A}$ in
the algebra of functions on $\mathcal{M}$ is less than or equal to the
norm of the element of $\mathcal{M}$.  In fact the norm of the
canonical homomorphism from $\mathcal{A}$ onto the algebra of
functions on $\mathcal{M}$ is equal to $1$, because the identity
element $e$ of $\mathcal{M}$ should have norm equal to $1$, and its
image is the constant function equal to $1$ on $\mathcal{M}$.

	Now suppose that our commutative algebra $\mathcal{A}$ is
equipped with an involution $*$, which is to say that there is a
conjugate-linear automorphism $x \mapsto x^*$ of the algebra such that
$(x^*)^* = x$ for all $x \in \mathcal{A}$.  Let $\phi$ be an element
of $\mathcal{M}$.  Using the involution $*$ we have that
$\overline{\phi(x^*)}$ is also a homomorphism from $\mathcal{A}$ onto
the complex numbers.  Depending on the situation, this may or may
not be the same as $\phi$.

	Let $x$ be an element of $\mathcal{A}$ which is in the kernel
of the canonical homomorphism from $\mathcal{A}$ onto the algebra of
functions on $\mathcal{M}$.  We can write $x$ as $x_1 + i \, x_2$
where $x_1$, $x_2$ are self-adjoint, i.e., $x_1^* = x_1$ and $x_2^* =
x_2$.  For each homomorphism $\phi$ from $\mathcal{A}$ onto the
complex numbers, we have that $\phi(x) = 0$, and therefore $\phi(x^*)
= 0$ as well, since $\overline{\phi(y^*)}$ is also a homomorphism from
$\mathcal{A}$ onto the complex numbers.  It follows that $\phi(x_1) =
\phi(x_2) = 0$, so that the kernel of the canonical homomorphism from
$\mathcal{A}$ onto the algebra of functions on $\mathcal{M}$ is
spanned by self-adjoint elements of $\mathcal{A}$.

	Suppose now that $V$ is a finite-dimensional complex vector
space of positive dimension equipped with a Hermitian inner product
$\langle v, w \rangle$.  Thus $v \mapsto \langle v, w \rangle$ is a
linear mapping from $V$ into the complex numbers for each $w \in V$,
$\langle w, v \rangle$ is equal to $\overline{\langle v, w \rangle}$
for all $v, w \in V$, and $\langle v, v \rangle$ is a nonnegative real
number for all $v \in V$ which is equal to $0$ if and only if $v = 0$.
If $T$ is a linear transformation on $V$, then there is a unique
linear transformation $T^*$ on $V$, called the adjoint of $T$, such
that
\begin{equation}
	\langle T(v), w \rangle = \langle v, T^*(w) \rangle
\end{equation}
for all $v, w \in V$.  The adjoint operation defines an involution
on the algebra $\mathcal{L}(V)$ of linear transformations on $V$.

	Let us suppose too that $\mathcal{A}$ is a commutative
subalgebra of $\mathcal{L}(V)$ which is invariant under the adjoint
operation, so that the adjoint defines an involution on $\mathcal{A}$.
Thus a self-adjoint element of $\mathcal{A}$ is a self-adjoint linear
transformation on $V$ in this situation.  The spectrum of a linear
operator consists of the eigenvalues of the operator, and the
eigenvalues of a self-adjoint operator are all real.  If $T$ is a
self-adjoint linear transformation on $V$ and $T^2 = 0$, then $T = 0$,
because $\langle T^2(v), v \rangle$ is equal to $\langle T(v), T(v)
\rangle$ in this case.

	Let $\mathcal{M}$ denote the set of homomorphisms from
$\mathcal{A}$ onto the complex numbers.  If $\phi \in \mathcal{M}$ and
$T$ is a self-adjoint element of $\mathcal{A}$, then $\phi(T)$ is a
real number.  This implies that $\phi(T^*) = \overline{\phi(T)}$ for
all $T \in \mathcal{A}$.  If $T$ is a self-adjoint element of
$\mathcal{A}$ which is nilpotent, then $T = 0$.  Therefore the
canonical homomorphism from $\mathcal{A}$ onto the algebra of
complex-valued functions on $\mathcal{M}$ is an isomorphism, and the
adjoint operation on $\mathcal{A}$ corresponds exactly to complex
conjugation of functions on $\mathcal{M}$.

	As an application, let $A$ be a finite commutative group, and
let $\mathcal{A}$ be the convolution algebra of complex-valued
functions on $\mathcal{A}$.  We can define an involution on
$\mathcal{A}$ by saying that the adjoint of a function $f(a)$ is given
by $\overline{f(-a)}$.  We can identify $\mathcal{A}$ with the algebra
of convolution operators on the space of functions on $A$, and this
definition of the involution is equivalent to taking the adjoint of
the convolution operator with respect to the usual inner product for
functions on $V$, where the inner product of functions $f_1$, $f_2$ is
the sum of $f_1(a) \, \overline{f_2(a)}$ over $a \in A$.  Therefore
the preceding discussion applies to $\mathcal{A}$, and since
homomorphisms from $\mathcal{A}$ onto the complex numbers correspond
exactly to homomorphisms of $A$ into the multiplicative group of
complex numbers with modulus $1$, it follows that the number of
homomorphisms of $A$ into the complex numbers with modulus $1$ is
equal to the number of elements of $A$.

	Suppose instead that $A$ is a finite group which may not be
commutative.  The corresponding convolution operator need not be
commutative either, but one can restrict oneself to the center of the
convolution algebra to get a commutative algebra.  Specifically, this
consists of the functions on $A$ which are constant on the conjugacy
classes of $A$, and the preceding discussion can also be applied to
the convolution algebra of these central functions on $A$.

\end{document}